


\parindent = 0 pt
\baselineskip = 16 pt
\parskip = \the\baselineskip

\font\AMSBoldBlackboard = msbm10

\def\RR{{\hbox{\AMSBoldBlackboard R}}}

\def\HH{{\hbox{\AMSBoldBlackboard H}}}

\settabs 12\columns
\rightline{15-JULY-1999}
\rightline{math.NT/9907103}
\vskip 1 true in
{\bf \centerline{ADDENDUM TO ``QUATERNIONIC GAMMA FUNCTIONS\dots''}}
\vskip 0.5 true in
\centerline{Jean-Fran\c{c}ois Burnol}
\par
\centerline{July 1999}
\vskip 0.5 true in
{\parskip = 0 pt
This note adds three annexes to my previous paper math/9904044\par
Annex 1. A sufficient condition for self-adjointness\par
Annex 2. Invariant closed operators on locally compact abelian groups\par
Annex 3. The trace of Connes for quaternions\par
This last item is a minor variation on the evaluation of Connes's trace (math/9811068), which is explained here in the setting of quaternions and can be applied also to any abelian local field.}
\par

\vfill
{\parskip = 0 pt
62 rue Albert Joly\par
F-78000 Versailles\par
France\par
\parskip = 12 pt
Address starting September 1st, 1999:\par
\baselineskip = 12 pt
Universit\'e de Nice-Sophia-Antipolis\par
\parskip = 0 pt
Laboratoire J.-A. Dieudonn\'e\par
Parc Valrose\par
F-06108 Nice C\'edex 02\par
France\par
}
The author thanks the ``Soci\'et\'e de Secours des Amis des Sciences'' for financial support during the time this work was completed.\par

\eject

The first two annexes can be read for the most part independently of my previous paper math/9904044, but the last one uses its notations with no further explanations.\par

{\bf Annex 1. A sufficient condition for self-adjointness\par}
{\bf Lemma:\quad} Let $L$ be a Hilbert space and $G$ a group of unitary operators on $L$. Let $A$ be the von Neumann algebra of bounded operators commuting with $G$. Let $M$ be a (possibly unbounded) operator with dense domain $D$. If\par
\+ \hfil (1) & $A$ is abelian\cr
\+ \hfil (2) & $M$ is symmetric\cr
\+ \hfil (3) & $M$ commutes with the elements of $G$\cr
then $M$ is essentially self-adjoint $\quad\bullet$\par

{\bf Note:\quad} Condition (3) is to be understood as follows: the dense domain $D$ is stable under $G$ and $M(g\cdot\varphi) = g\cdot(M(\varphi))$ holds for $\varphi \in D$ and $g \in G$. \par

{\bf Proof:\quad} We first replace $(M,D)$ by its double-adjoint so that we can assume that $(M,D)$ is closed (it is easy to check that conditions (2) and (3) remain valid). The problem is to show that it is self-adjoint. Let $K$ be the range of the operator $M + i$. It is a closed subspace of $L$ (as $\|(M + i)(\varphi)\|^2 = \|M(\varphi)\|^2 + \|\varphi\|^2$, and $M$ is closed). Let $R$ be the bounded operator onto $D$ which is orthogonal projection onto $K$ followed with the inverse of $M + i$. It belongs to $A$, hence commutes with its adjoint $R^*$. Any vector $\psi$ in the kernel of $R$ is then in the kernel of $R^*$ (as $<R^*\psi|R^*\psi>\ =\ <\psi|R\,R^*\psi>\ =\ 0$). So $\psi$ belongs to the orthogonal complement to the range of $R$, that is $\psi = 0$ as the range of $R$ is $D$. So $K = L$ and in the same manner $(M - i)(D) = L$. By the basic criterion for self-adjointness ([ReSi80]), $M$ is self-adjoint$\quad\bullet$\par

{\bf Example:\quad} The argument leading to the proof of Theorem I shows that the von Neumann algebra of bounded operators commuting simultaneously with the left and right actions of the multiplicative quaternions on $L^2(\HH, dx)$ is abelian. So the lemma can be applied to the operator $\log(|x|) + \log(|y|)$, with initial domain the space of Schwartz functions (for example).\par

{\bf Annex 2. Invariant closed operators on locally compact abelian groups\par}
Let $G$ be a locally compact, Hausdorff, topological abelian group. We refer to [Ru62] for the basics of harmonic analysis on $G$. In particular we have a Haar measure $dx$ and a Hilbert space $L = L^2(G, dx)$. We also assume that $G$ is a countable union of subsets with finite measure, so that there exists $\psi \in L$ with the property $\psi(x)\neq 0\ (a.e)$. We will use $\psi(x)$ below.\par

The dual group $\widehat{G}$ of unitary characters acts on $L$ by multiplication operators (i.e. $(\chi\cdot\varphi)(x) = \chi(x)\varphi(x)$), and we are interested in operators (possibly unbounded) that commute with this action (in the sense explained in Annex 1). The lemma we are aiming at is certainly extremely well-known, but I could not find an appropriate reference.\par

Let $a(x)$ be a measurable function on $G$, not necessarily bounded. Let $D_a$ be the domain of elements $\varphi$ of $L$ such that $a(x)\varphi(x)$ also belongs to $L$. And let $M_a$ be the operator with domain $D_a$ acting through multiplication with $a$. We note that $D_a$ is stable under $\widehat{G}$ and that $(M_a, D_a)$ commutes with $\widehat{G}$. Furthermore $D_a$ is dense: first it contains $$\psi(x)\over\sqrt{1 + |a(x)|^2}$$
Hence if $f$ is orthogonal to $D_a$ then the function defined as $$\overline{f(x)}\psi(x)\over\sqrt{1 + |a(x)|^2}$$ belongs to $L^1(G,dx)$ and has a vanishing ``Fourier transform'', hence $f = 0$. It is also clear that $(M_a, D_a)$ is a closed operator (as a sequence converging in the mean has a pointwise a{.}e{.} convergent subsequence). One more remark is that if for another function $b(x)$ the operator $(M_b, D_b)$ extends $(M_a,D_a)$, then in fact $a = b\ (a.e.)$ and $(M_b, D_b) = (M_a,D_a)$. Finally we note that the adjoint of $(M_a,D_a)$ is $(M_{\overline{a}},D_{\overline{a}})$. This can be seen as follows: Let $f(x)$ be in the domain of the adjoint. There exists then an element $\theta$ of $L$ such that for any $\varphi \in D_a$ the equality $$\int f(x)\overline{a(x)\varphi(x)}\,dx = \int \theta(x)\overline{\varphi(x)}\,dx$$ holds. This implies that the two functions of $L^1(G, dx)$ given as
$${f(x)\overline{a(x)\psi(x)}\over\sqrt{1 + |a(x)|^2}}\hbox{\quad and\quad}{\theta(x)\overline{\psi(x)}\over\sqrt{1 + |a(x)|^2}}$$
have the same Fourier transform, hence are equal almost everywhere. So $f \in D_{\overline{a}}$ and $(M_a)^*(f) = (M_{\overline{a}})(f)$\par

Our ``certainly extremely well-known'' lemma now reads as:\par

{\bf Lemma:\quad} Let $(M, D)$ be a closed operator commuting with the action of $\widehat{G}$. Then $(M, D) = (M_a, D_a)$ for a (unique) measurable function $a(x)$ (two such functions being identified if they are equal almost everywhere)$\bullet$

{\bf Note:\quad} for a bounded $M$ and $G = \RR$, this is proven in the classical reference [StWe73], as a special case of a more general statement in $L^p$ spaces. Their proof uses distributions.\par

{\bf Proof:\quad} Let us first assume that $M$ is bounded. We then pick as above some $\psi(x)$ non-vanishing almost everywhere and define $a(x) = {M(\psi)(x)\over\psi(x)}$. Let us consider the domain $D$ consisting of all finite linear combinations of $\chi(x)\psi(x)$, $\chi \in \widehat{G}$. It is dense by the argument about unicity of Fourier transform in $L^1$ we have used many times. Then $(M, D) \subset (M_a, D_a)$, hence $(M_a, D_a)$ is also an extension of the closure of $(M, D)$ which is $(M, L)$. But this means that $D_a = L$ and that $M = M_a$ (as a further token we note that necessarily $a$ is essentially bounded)\par

The next case is when $M$ is assumed to be self-adjoint. Its resolvents $R_1 = (M - i)^{-1}$ and $R_2 = (M + i)^{-1} = R_1^*$ are bounded and commute with $\widehat{G}$. Hence they correspond to multiplicators $r_1(x)$ and $r_2(x) = \overline{r_1(x)}$. The kernel of $R_1$ is orthogonal to the range of $R_2$ which is all of $D$, so that it is reduced to $\{0\}$. So $r_1(x) \neq 0\ (a.e.)$. Let $f \in D$ and $g = M(f)$.
$$R_1(M(f) - i\,f) = f$$
$$r_1(x)(g(x) - i\,f(x)) = f(x)$$
$$g(x) = {1 + i\,r_1(x)\over r_1(x)}\cdot f(x)$$
So $(M_a, D_a)$ is an extension of $(M, D)$ for $a(x) = {1 + i\,r_1(x)\over r_1(x)}$. Switching to the adjoints we deduce that $(M, D)$ is an extension of $(M_{\overline{a}},D_{\overline{a}})$. So all three are equal (and $a$ is real-valued).\par

For the general case we use the theorem of polar decomposition ([ReSi80]). There exists a non-negative self-adjoint operator $|M|$ with the same domain as $M$ and a partial isometry $U$ such that $M = U|M|$. Further conditions are satisfied which make $|M|$ and $U$ unique: so they also commute with $\widehat{G}$. It follows from what was proven previously that $(M, D) \subset (M_a, D_a)$ for an appropriate $a$ (the product of the multiplicators associated to the self-adjoint $|M|$ and the bounded $U$). The adjoint $(M^*, D^*)$ also has a dense domain and commutes with $\widehat{G}$, so in the same manner $(M^*, D^*) \subset (M_b, D_b)$ for an appropriate $b$. The inclusion $(M_{\overline{a}},D_{\overline{a}}) \subset (M^*, D^*) \subset (M_b, D_b)$ implies $b = \overline{a}$ and $(M_a, D_a) = (M, D)^{**}$. But the double-adjoint coincides with the closed operator $(M,D) \quad\bullet$\par

{\bf Note:\quad} it follows that a closed symmetric operator commuting with $\widehat{G}$ is self-adjoint, and that a symmetric operator with dense domain that commutes with $\widehat{G}$ is essentially self-adjoint (alternatively Annex 1 applies here).

{\bf Annex 3. The trace of Connes for quaternions\par}

Let $f(g)$ be a smooth function with compact support on $\HH^\times$. Let $U_f$ be the bounded operator $\int f(g) L_2(g)\,d^*g$ on $L = L^2(\HH, dx)$. So
$$U_f: \varphi(x) \mapsto \int f(g){1\over\sqrt{|g|}}\varphi(g^{-1}x)\,d^*g$$
The composition $U_f\,{\cal F}$ of $U_f$ and the Fourier Transform ${\cal F}$ acts as
$$\eqalign{
\varphi(x) &\mapsto \int\int f(g){1\over\sqrt{|g|}}\lambda(-g^{-1}xy)\varphi(y)\,d^*g\,dy \cr
&= {1\over\sqrt{2\pi^2}}\int \left( \int f({1\over g}){1\over\sqrt{2\pi^2|g|}}\lambda(-gxy)\,dg\right) \varphi(y)\,dy \cr
&= {1\over\sqrt{2\pi^2}}\int \widetilde{I(f)_a}(xy)\varphi(y)\,dy\cr
}$$
In this last equation $\widetilde{I(f)_a}$ is the additive Fourier Transform of the additive representative of $I(f)(g)$ (i.e. $f({1\over g}){1\over\sqrt{2\pi^2|g|}}$).\par

Following Connes [Co98], our goal is to compute the trace $T(\Lambda)$ of the operator $\widetilde{P_\Lambda}P_\Lambda\,U_f$, where $\widetilde{P_\Lambda} = {\cal F}P_\Lambda{\cal F}^{-1}$ and $P_\Lambda$ is the cut-off projection to functions with support in $|x| \leq \Lambda$. Connes evaluated this trace not only on a local field, but also in a situation involving his crucial quotient construction (when more than one place are considered). I will now explain in the set-up of quaternions a method which would apply almost word for word to the simple case of a local field (when the quotient construction is not involved). It is a minor variation on Connes's method, avoids the use of distributional kernels and symbols, and checks explicitely that the operator considered is trace-class. Our reference for trace-class operators is [GoKr69].\par

We recall that if $A$ is trace-class then for any bounded $B$, $AB$ and $BA$ are trace-class and have the same trace. Also if $K_1$ and $K_2$ are two Hilbert-Schmidt operators given for example as $L^2-$kernels $k_1(x,y)$ and $k_2(x,y)$ on a measure space $(X, dx)$ then $A = K_1^*\, K_2$ is trace-class and its trace is the Hilbert-Schmidt scalar products of $K_1$ and $K_2$:
$$Tr(K_1^*\, K_2) = \int\int \overline{k_1(x,y)}\,k_2(x,y)\ dxdy$$

The operator $P_\Lambda {\cal F}^{-1} P_\Lambda$ is an operator with as kernel a smooth function restricted to a finite box (precisely it is $\lambda(xy),\ |x|,|y|\leq\Lambda$). Such an operator is trace class, as is well-known (one classical line of reasoning is as follows: taking a smooth function $\rho(x)$ with compact support, identically $1$ on $|x|\leq\Lambda$, and $Q_\rho$ the multiplication operator with $\rho$, one has $P_\Lambda {\cal F}^{-1} P_\Lambda = P_\Lambda Q_\rho {\cal F}^{-1} Q_\rho P_\Lambda$, so that it is enough to prove that $Q_\rho {\cal F}^{-1} Q_\rho$ is trace-class. This operator has a smooth kernel with compact support, so we can put the system in a box, and reduce to an operator $K$ with smooth kernel on a torus. Then $K = (1 + \Delta)^{-n}(1 + \Delta)^{n}K$ with $\Delta$ the positive Laplacian. For $n$ large enough, $(1 + \Delta)^{-n}$ is trace-class, while $(1 + \Delta)^{n}K$ is at any rate bounded.)\par

So Connes's operator $\widetilde{P_\Lambda}P_\Lambda\,U_f = {\cal F}\cdot P_\Lambda {\cal F}^{-1} P_\Lambda\cdot U_f$ is indeed trace class and 
$$T(\Lambda) = Tr(\widetilde{P_\Lambda}P_\Lambda\,U_f) = Tr(P_\Lambda {\cal F}^{-1} P_\Lambda\cdot U_f{\cal F}) = Tr(P_\Lambda {\cal F}^{-1} P_\Lambda\cdot P_\Lambda U_f{\cal F} P_\Lambda)$$
can be computed as a Hilbert-Schmidt scalar product:
$$T(\Lambda) = {1\over\sqrt{2\pi^2}}\int\int_{|x|, |y| \leq\Lambda}\lambda(xy)\widetilde{I(f)_a}(xy)\,dxdy$$
using the change of variable $(x,y) \mapsto (z= xy, y)$
$$T(\Lambda) = \sqrt{2\pi^2} \int_{|z|\leq\Lambda^2} \lambda(z) \widetilde{I(f)_a}(z)\left(\int_{{|z|\over \Lambda}\leq|y|\leq\Lambda} {dy\over 2\pi^2 |y|}\right)\,dz$$
$$T(\Lambda) = \sqrt{2\pi^2} \int_{|z|\leq\Lambda^2} (2\log(\Lambda) - \log(|z|))\lambda(z) \widetilde{I(f)_a}(z)\,dz\leqno(C)$$
With the notation $B = \log(|z|) = {\cal F}^{-1}\log(|x|){\cal F}$ and $(2\log(\Lambda) - B)_+ = Max(2\log(\Lambda) - B, 0)$, and recalling that $\sqrt{2\pi^2}$ is involved in switching from the additive to multiplicative picture, we can conclude:

{\bf Lemma:\quad}$\widetilde{P_\Lambda}P_\Lambda\,U_f$ is trace-class and
$$\eqalign{
Tr(\widetilde{P_\Lambda}P_\Lambda\,U_f) &= (2\log(\Lambda) - B)_+(I(f))(1)\cr
&= 2\log(\Lambda)f(1) - H(f)(1) + o(1)\cr
}$$

For the last line we used that $B(I(f))(1) = H(I(f))(1) = H(f)(1)$ as $H = \log(|x|) + \log(|z|)$ commutes with the Inversion $I$. The error (for $\Lambda \rightarrow\infty, \Lambda > 1$) is $o(1)$ as it is bounded above by $$\sqrt{2\pi^2} \int_{|z|\geq\Lambda^2} (\log(|z|)) \left|\widetilde{I(f)_a}(z)\right|\,dz$$ and $\widetilde{I(f)_a}(z)$ is a Schwartz function.\par

It is interesting to note that in Connes's paper (in an abelian situation, involving either one place or more than one in the context of his quotient) the computation also goes through an intermediate stage essentially identical with $(C)$ and that the identification of the constant term $H(f)(1)$ with a quantity related to the explicit formula of number theory then requires a further discussion, going through an intermediate stage (Weil's form of the local terms of the explicit formula). The main result of this paper (math/9904044) and of the previous one [Bu99] is the direct connection between $H$ and the logarithmic derivatives of the Tate Gamma functions (which are involved in the analytic expression of the explicit formula), and allows to avoid this ``detour''. But Weil's form seems to play a crucial motivational r\^ole in Connes's paper through the connection with a fixed point formula.\par

{{\bf REFERENCES}\par
\font\smallRoman = cmr10
\smallRoman
\font\smallBold = cmbx10
\font\smallSlanted = cmsl10

{\smallBold [Bu99] J.F. Burnol}, {\smallSlanted ``The Explicit Formula and the conductor operator''}, math/9902080 (February 1999).\par

{\smallBold [Co98] A. Connes}, {\smallSlanted ``Trace formula in non-commutative Geometry and the zeros of the Riemann zeta function''}, math/9811068 (November 1998).\par

{\smallBold [GoKr69] I. Ts. Gohberg, M. G.Krein}, {\smallSlanted ``Introduction to the theory of linear non-self-adjoint operators''}, American Mathematical Society, (1969), (original in russian 1965).\par

{\smallBold [ReSi80] M. Reed, B. Simon}, {\smallSlanted ``Methods of modern mathematical physics, Vol.1: Functional Analysis''}, revised and enlarged edition, Academic Press (1980).\par

{\smallBold [Ru62] W. Rudin}, {\smallSlanted ``Fourier analysis on groups''}, Interscience Publishers, (1962).\par

{\smallBold [StWe73] E. Stein, G. Weiss}, {\smallSlanted ``Introduction to Fourier Analysis on Euclidean Spaces''}, Princeton University Press, (1973).\par

Jean-Fran\c{c}ois Burnol, 62 rue Albert Joly, F-78000 Versailles, France, July 1999.\par
}
\vfill
\eject
\bye